\documentclass[a4paper]{article}
\input{INF.STY}
\theoremstyle{plain}
\begin{document}
 \author{ Е. В. Щепин
 \thanks{ Работа написана в рамках программы РАН
 "Алгебраические и комбинаторные методы математической кибернетики".}}

 \title{О фрактальных кривых Пеано}
\date{\it Посвящаяется памяти Людмилы Всеволодовны Келдыш}

\maketitle

\paragraph{Введение.}
Под \emph{кривой Пеано} подразумевается любое непрерывное
отображение числового отрезка на плоский квадрат.  В \cite{TULA87}
обоснована целесообразность применения кривых Пеано для построения
развертки телевизора с электронно-лучевой трубкой.\footnote{ В
реальности всегда, конечно, речь идет не о непрерывной кривой
Пеано, а о каком то ее дискретном приближении.} В пеановских
телевизорах помехи сворачиваются, по сравнению с обычной строчной
разверткой, что должно приводить к большей устойчивости
изображения. Развертка кадра  представляет собой кривую ---
траекторию светящейся точки на экране монитора. Число светящихся
точек на экране конечно и электронный луч по очереди освещает эти
точки. Поэтому развертка представляет собой цифровую кривую ---
отображение цифрового (т.е. состоящего из конечного числа точек
--- пикселей) отрезка на цифровой прямоугольник --- экран
телевизора. Помеху представим себе как кратковременное нарушение
передачи сигнала. Если помеха произошла в некоторый интервал
времени, то на экране нарушится изображение на развертке этого
интервала, то есть в той области экрана, которую электронный луч
проходит за этот временной интервал. Площадь помехи (количество
пикселей экрана в ней содержащихся) определяется величиной
временного интервала помехи и зависит лишь от частоты кадров. Она
будет одинакова как для построчной, так и для пеановской
разверток, поскольку и та и другая за одинаковое время (время
между двумя кадрами) проходят равное число пикселей (количество
всех пикселей экрана). Но форма помехи будет существенно
различаться. Если для построчной развертки помеха вытягивается
вдоль строк через весь экран, то в случае пеановской развертки она
сворачивается и может иметь весьма небольшой диаметр.  Среди
различных разверток предпочтение следует отдать той, которая в
наибольшей степени сворачивает помехи. Вытянутость помехи может
быть охарактеризована отношением квадрата ее диаметра к ее
площади.  Если кривая Пеано заметает площадь с постоянной
скоростью (а таковы все фрактальные кривые Пеано), то способность
сворачивания помех кривой $p(t)$ характеризуется максимумом
\emph{квадратно-линейного отношения}
$\frac{|p(t)-p(t')|^2}{|t-t'|}$. Например, классическая кривая
Пеано-Гильберта, как показал К.Бауман \cite{B}, имеет максимальное
квадратно-линейное отношение равное шести. Открытым является
следующий вопрос.
\begin{prob}Каково наименьшее число $\kappa$, для которого
существует непрерывное отображение единичного отрезка на единичный
квадрат с максимальным квадратно-линейным отношением равным
$\kappa$?
\end{prob}
Из теоремы Арцела вытекает компактность множества отображений
единичного отрезка в единичный квадрат, имеющих конечное
квадратно-ли\-ней\-ное отношение. Поэтому число $\kappa$
существует. Это число называем \emph{константой Пеано}.  Из
результатов статьи \cite{9} вытекает, что $\kappa\le 5\frac23$. В
настоящей работе доказано, что $\kappa\ge 3\frac12$. А основным
результатом настоящей статьи является доказательство того факта,
что для широкого класса фрактальных пеановских кривых максимальное
квадратно-линейное отношение не может быть менее пяти.

%

\paragraph{Нижняя оценка константы Пеано.}

Под \emph{кривой}, мы всюду ниже понимаем непрерывное отображение
в плоскость числового интервала, который мы будем
интерпретировать, как  временной и элементы которого будут
называться \emph{моментами} времени.

 \emph{Единичной} пеановской кривой
называется  кривая,  отображающая единичный отрезок на единичный
квадрат.

\begin{tr} Квадратно-линейное отношение единичной пеановской кривой,
не может быть меньше чем $3.5$
\end{tr}
\begin{proof}
Обозначим через $A,B,C,D$ вершины квадрата, а через $a,b,c,d$
моменты времени, в которые их проходит кривая. Предположим, что
$a<b<c<d$. Теперь, во первых, заметим, что пары вершин $A,B$,
$B,C$, $C,D$ являются соседними вершинами квадрата. Действительно,
если квадратно-линейное отношение не превосходит $3.5$, то
расстояния между прообразами соседних вершин не может быть меньше
чем $2/7$, а между прообразами противоположных вершин не может
быть меньше чем $4/7$. Поэтому сумма длин
$|b-a|+|c-b|+|d-c|=|d-a|\le1$ может быть меньше единицы, только
если все слагаемые  меньше чем $4/7$. Далее заметим, что прообраз
отрезка $AD$ не пересекает отрезка $[b,c]$, иначе длина этого
отрезка была бы больше чем $4/7$. Поскольку прообразы концов
отрезка $AD$ лежат по разные стороны отрезка $[b,c]$, постольку
найдется точка отрезка $AD$, на расстоянии $x$ от $A$, в любой
окрестности которой имеются точки, лежащие по разные стороны от
$[b,c]$. Тогда, если квадратно-линейное отношение не превосходит
$3.5$, то получим два неравенства на $x$.
\begin{gather}\label{2in}
(1+x^2)\le 3.5 b\\1+(1-x)^2\le 3.5(1-c).
\end{gather}
 Складывая эти
неравенства, получаем $2+x^2+(1-x)^2\le 3.5-3.5(c-b)$. Поскольку
левая часть этого неравенства достигает наименьшего значения
равного $2.5$ при $x=1/2$, а правая часть никак не может быть
больше $2.5$, ибо $c-b\ge2/7$, постольку последнее неравенство
может быть справедливым только при $x=1/2$ и $c-b=2/7$. Тогда из
первого неравенства \eqref{2in} следует $b\ge 5/14$, а из второго
--- $c\le 9/14$. Откуда  $c-b\le 2/7$, а значит и квадратно-линейное отношение
 отображения, переводящего $b,c\to B,C$, не меньше чем
$3.5$.
\end{proof}

\paragraph{Изометричные и подобные кривые.}
 \emph{ Изометрией} между
кривыми $f\:[a,b]\to Q_1$ и $g\:[c,d]\to Q_2$ называется пара
изометрий $h,h'$, где $h$ --- изометрия прямой, а $h'$
--- изометрия плоскости, такая что $h'f=gh$.


 Две кривые $f$, $g$ называются
 \emph{подобными} c коэффициентом подобия $k$, если кривая
 $\sqrt kf(t/k)$ изометрична кривой $g(t)$. Например, если $k=2$,
 то $g(t)$ определена на отрезке длины вдвое большем, чем $f(t)$ и имеет образ
 площади вдвое большей чем  $f(t)$.

\paragraph{Фрактальные периоды.} Отрезок,
 содержащийся в области определения кривой $f(t)$,
называется \emph{фрактальным периодом}
 этой кривой, если ограничение
кривой на этот отрезок подобно всей кривой. В дальнейшем
\emph{период} по умолчанию означает \emph{фрактальный период.}

 Период называется
\emph{максимальным}, если он не содержится ни в каком другом
периоде. \emph{Порядком периода} называется длина максимальной
растущей цепочки периодов, начинающейся с данного. Ограничение
кривой на ее период называется \emph{фракцией}. Соответственно
\emph{максимальной фракцией} называется ограничение кривой на
максимальный период. И порядок фракции определяется как порядок ее
периода.

\begin{lem}
\label{isolat} Всякое замкнутое несвязное подмножество $S$ отрезка
$[a,b]$ содержит как точку изолированную слева так и точку
изолированную справа.
\end{lem}
\begin{proof}
Так как $S$ несвязно, найдется точка $c\in[a,b]$, слева и справа
от которой имеются точки $S$. Тогда $\max \{s\in S\mid s<c\}$
является точкой $S$ изолированной справа, а $\min \{s\in S\mid
s>c\}$ --- изолированной слева.
\end{proof}

\begin{lem}
Всякая растущая цепочка периодов обрывается. Всякая убывающая
последовательность периодов сходится к нулю.
\end{lem}
\begin{proof}
Рассмотрим множество коэффициентов подобия кривой $f$ своим
фракциям. Это множество во-первых замкнуто. Действительно, если
$k_1,k_2,\dots $ последовательность коэффициентов подобия и
$T_1\subset T_2\subset\dots T_n\subset$ соответствующая им
последовательность периодов      кривой $f$. Пусть $h_i,h_i'$
последовательность подобий, переводящая $f$ в ограничение $f$ на
$T_i$. Тогда, переходя, если нужно к подпоследовательностям, можем
считать, что $h_i$ сходятся к $h$, и $h_i'$ сходятся к $h'$, а
периоды сходятся к $T$. В таком случае ограничение $f$ на $T$
подобно $f$ посредством пары подобий $h,h'$, и коэффициент
подобия, очевидно, равен пределу $k_i$.

Если мы предположим, что имеется бесконечная возрастающая
последовательность периодов, то в множестве коэффициентов подобия
все точки будут неизолированными слева. А замкнутое множество с
таким свойством связно в силу леммы \ref{isolat}. В этом случае
кривая обладает непрерывным подобием, что означает, что кривая
изометрична линейному отображению.

Аналогично, если имеется убывающая последовательность
коэффициентов подобия, не сходящаяся к нулю, то все элементы
множества коэффициентов неизолированы справа, что тоже влечет
связность этого множества и линейность кривой.
\end{proof}

\begin{corr} Всякий фрактальный период имеет конечный порядок.
В частности, всякий фрактальный период  содержится в максимальном.
\end{corr}

\paragraph{Правильные разбиения.}
Последовательность равных отрезков $I_1,\dots I_k$ называется
\emph{правильным разбиением отрезка} $I$, если $I=\cup I_k$ и
различные отрезки последовательности  имеют дизъюнктные
внутренности.

\emph{Правильным разбиением} квадрата $Q=I\times I$ называется
совокупность произведений $\{I_i\times I_j\}_{i,j\le k}$, где
$I_1,\dots I_k$ является правильным разбиением отрезка $I$.

\begin{lem}\label{regpart}
Если  квадрат $Q$ представлен в виде объединения равных квадратов
$\{Q_i\}_{i\le n}$, суммарная площадь которых равна площади $Q$,
то $\{Q_i\}_{i\le n}$ является правильным разбиением $Q$.
\end{lem}
\begin{proof}
Предоставляется читателю.
\end{proof}

\paragraph{Фрактальные пеановские кривые.}
Кривая $f(t)$ называется \emph{фрактальной}, если каждая точка из
области ее определения принадлежит некоторому фрактальному
периоду.

Фрактальная  кривая называется \emph{правильной пеановской}, если
ее образ является квадратом и максимальные периоды образуют
правильное разбиение области определения.

\begin{lem}
Образы максимальных фракций правильной пеановской кривой $f\:
[a,b]\to Q$ образуют правильное разбиение квадрата $Q$. В
частности количество максимальных фракций является полным
квадратом.
\end{lem}
\begin{proof}
Это немедленно вытекает из Леммы \ref{regpart}.
\end{proof}

  Мы будем говорить, что
правильная пеановская кривая имеет \emph{фрактальный род} $k$,
если она имеет $k$ максимальных периодов.

 Классическая кривая
Пеано-Гильберта является правильной пеановской кривой фрактального
рода $4$. А оригинальная кривая Пеано  имеет фрактальный род $9$.

\paragraph{Конечность квадратно-линейного отношения.}
\emph{Стыком} фрактальной кривой называется кривая, полученная
ограничением на пару соседних фрактальных периодов (любого
порядка).
\begin{lem}\label{styk-finite} Для любой правильной
пеановской кривой имеется лишь конечное число попарно
неизометричных стыков.
\end{lem}
\begin{proof} Предоставляется читателю.
\end{proof}
Если $f\:[a,b]\to \C$ является стыком кривой фрактального рода
$k$, то определим \emph{квадратно-линейное отношение стыка} $f$,
как максимум квадратно-линейных  отношений для пар $x,y\in[a,b]$,
таких что $|x-y|\ge \frac{ |b-a|}{2k}$.
\begin{lem}\label{styk-max} Квадратно-линейное отношение правильной пеановской
кривой равно максимуму из квадратно-линейных отношений ее стыков.
\end{lem}
\begin{proof}
Рассмотрим произвольную пару моментов $t,t'$. Пусть $k$ ---
наименьшее число, для которого $t,t'$ не принадлежат стыку
фрактальных периодов $n$-го порядка. Тогда $t-t'$ больше, чем
длина периода $n$-го порядка, поэтому пара $t,t'$ принадлежит
стыку $(n-1)$-го порядка и участвует в определении его
квадратно-линейного отношения.
\end{proof}

\begin{tr}\label{finitness} Для любой правильной фрактальной
пеановской кривой квадратно-линейное отношение ограничено сверху.
\end{tr}
\begin{proof}
Немедленно вытекает из лемм \ref{styk-finite} и \ref{styk-max}.
\end{proof}

\paragraph{Диагональные и односторонные кривые.}В соответствии
 с доказанной ниже теоремой \ref{begin-end} правильные пеановские
кривые делятся на два типа: \emph{диагональные}, у которых начало
и конец лежат на диагонали квадрата и \emph{односторонние}, у
которых начало и конец лежат на одной стороне.

\begin{lem}
Если $k>2$, то невозможно перенумеровать все квадраты
$\{Q_i\}_{i\le k^2}$ правильного разбиения квадрата $Q$ так, чтобы
квадраты $Q_i$ и $Q_{i+1}$ при $i=1,2,\dots,k^2-1$ пересекались по
стороне, а квадраты $Q_i$ и $Q_{i+2}$ при $i=1,2,\dots,k^2-2$
пересекались по вершине.
\end{lem}
\begin{proof}
Предположим, что такая нумерация существует. Рассмотрим тот из
угловых квадратиков $Q$, скажем $Q_i$, для которого $k^2-1>i>2$.
Тогда $Q_{i-2},Q_{i-1},Q_{i+1}, Q_{i+2}$ должны иметь непустое
пересечение с $Q_i$. Но для углового квадратика есть только три
соседних квадратика, которые его пересекают.
\end{proof}
\begin{tr}\label{begin-end}
Для любой правильной пеановской кривой $f\:[a,b]\to Q$ ее начало
$f(a)$ и конец $f(b)$ являются различными вершинами квадрата $Q$.
\end{tr}

\begin{proof}
Предположим, что начало кривой  $p(a)$ не является вершиной
квадрата $Q$. Пусть $\{Q_i\}_{i=1}^{k^2}$ является правильным
разбиением на образы максимальных фракций кривой $p(t)$,
занумерованное в порядке их прохождения кривой.  Положим
$b_1=a+(b-a)/k^2$. Тогда $p(b_1)$ принадлежит пересечению $Q_1\cap
Q_2$ и не совпадает с вершинами $Q_i$, потому что $b_1$ является
началом второй фракции. С другой стороны $b_1$ является концом
первой фракции, поэтому получаем, что и $p(b)$ не является
вершиной $Q$. Рассмотрим различные варианты расположения $p(a)$ и
$p(b)$ на сторонах $Q$.

1) $p(a)$ и $p(b)$ принадлежат одной и той же стороне. Это
невозможно потому что в таком случае $Q_2$ должен пересекать $Q_1$
и $Q_3$ по одной и той же стороне.

2) $p(a)$ и $p(b)$ лежат на противоположных сторонах $Q$. В этом
случае  цепочка кубиков $\{Q_i\}$ вытягивается в линию и не может
образовать квадрат.

3)  $p(a)$ и $p(b)$ лежат на соседних сторонах квадрата. В этом
случае пересечение $Q_i$ с $Q_{i+2}$ непусто и можно сослаться на
лемму.

\end{proof}

\paragraph{Угловые моменты.}
Момент времени $\tau$ называется \emph{угловым порядка $k$} для
правильной пеановской кривой $p(t)$, если $p(\tau)$ является
вершиной образа фракции $k$-го порядка.

 Для каждого фрактального периода  его начало и конец являются
угловыми моментами  в силу теоремы \ref{begin-end}. В силу леммы
\ref{vertex} каждый фрактальный период порядка $k$, помимо начала
и конца содержит ровно два угловых момента $k$-го порядка, которые
мы будем называть \emph{внутренними}, а их образы называть
внутренними вершинами образов фракций: первой и второй в
соответствии с порядком их прохождения кривой.

\begin{lem}
\label{vertex} Если $f$ --- фрактальная пеановская кривая, образом
которой является квадрат $Q$, то прообраз любой из вершин квадрата
$Q$ состоит ровно из одной точки.
\end{lem}
\begin{proof}
Для любой вершины и любого $k$ существует единственная фракция
$k$-го порядка, содержащая эту вершину. Поэтому прообраз вершины
целиком лежит во фрактальном периоде $k$-го порядка. Так как это
верно при любом $k$, то прообраз одноточечен.
\end{proof}

\paragraph{Пятерные кривые.} Правильную пеановскую кривую мы будем
называть \emph{пятерной}, если она подобна единичной пеановской
кривой и имеет квадратно-линейное отношение менее пяти. Нашей
конечной целью является доказательство несуществования пятерных
кривых.

\begin{lem}
\label{2-diag} Если кривая $f\:[0,1]\to Q$ пятерная, то вторая,
проходимая ей  вершина $Q$, лежит на одной стороне с первой
(точкой входа), а третья ---  на одной стороне с четвертой (точкой
выхода).
\end{lem}
\begin{proof}
 Пусть $ABCD$
--- вершины (единичного) квадрата $Q$, указанные в порядке их прохождения
кривой $f$, а $0=a,b,c,d=1$ соответственно обозначают моменты
прохождения этих вершин кривой $f$. Так что $f(a)=A$
--- начало, а $f(d)=D$
--- конец кривой. Если $B=f(b)$ --- первая внутренняя вершина ---
лежит на одной стороне с $A$, то $C=f(c)$ не остается ничего
иного, как лежать на одной стороне с $D$. То есть в этом случае
заключение леммы справедливо. Если же $B$ лежит с $A$ на одной
диагонали, то пятерность кривой дает неравенство $b-a>2/5$. Далее
в этом случае $C$ лежит на одной диагонали с $D$, откуда
$d-c>2/5$. Наконец, расстояние между $B$ и $C$ единица, откуда
$c-b>1/5$. Складывая вместе все полученные неравенства, получаем
$d-a>1$ вопреки условию $d-a=1$.
\end{proof}

\begin{lem}\label{str-cont}
Если  кривая $f$ пятерная, то образы любой пары соседних периодов
одного порядка являются квадратами с общей стороной.
\end{lem}
\begin{proof}
Предположим, что образы некоторой пары соседних периодов одного
порядка кривой $f(t)$ пересекаются по вершине и докажем, что в
таком случае $f$ не является пятерной.

Переходя от кривой $f$ к подобной кривой $\phi(t)=\sqrt kf(kt)$,
мы, не меняя квадратно-линейного отношения, при соответствующем
подборе $k$, можем добиться, чтобы образы пары соседних периодов
кривой $\phi$ представляли собой единичные квадраты $Q_i$ и
$Q_{i+1}$, пересекающиеся по вершине, а соответствующие им периоды
имели единичную длину.

Пусть $t_0,t_1,t_2,t_3,t_4,t_5,t_6$ угловые моменты
последовательного прохождения вершин этой пары квадратов кривой
$\phi$. То есть $t_0$ --- начало $i$-го периода (момент вхождения
в $Q_i$), $t_6$
--- конец $i+1$-го (момент выхода из $Q_{i+1}$), а
$t_3$, момент перехода из $Q_i$ в $Q_{i+1}$, то есть конец $i$-го
периода, совпадающий с началом $i+1$-го.

Пусть кривая $f$ имеет диагональный тип. Тогда
$|\phi(t_1)-\phi(t_0)|=1$ и $|\phi(t_2)-\phi(t_1)|^2=2$ и, так как
квадратно-линейное отношение  $\phi$ не превышает пяти, то
$t_1-t_0\ge\frac15$ и $t_1-t_0\ge\frac25$. Откуда получаем
$t_3-t_2\le\frac25$. Аналогично доказывается, что $t_4-t_3\le
\frac25$. Получаем $t_4-t_2\le\frac45$. А так как $|\phi
(t_4)-\phi(t_2)|=2$, то квадратно-линейное отношение кривой $\phi$
на паре $t_4,t_2$ больше либо равно пяти.

Если кривая односторонняя и имеет квадратно-линейное отношение
$\le 5$, то $t_1-t_0\ge \frac15$ и $t_6-t_5\ge \frac15$, откуда
$t_5-t_1\le \frac 85$. А так как $|f(t_5)-f(t_1)|^2=8$, получаем,
что квадратно-линейное отношение кривой $f$ на паре $t_1,t_5$
больше либо равно пяти.
\end{proof}

\paragraph{Ускорение и замедление.}
 Кривую
будем называть \emph{ускоряющейся}, если промежуток времени от
начала до первого внутреннего углового момента превосходит
промежуток времени от второго внутреннего углового момента до
конца и называем \emph{замедляющейся} в случае  противоположного
неравенства между этими промежутками.

\emph{Смещением} кривой $f$ от момента $t$ до момента $t'$
называем вектор $f(t')-f(t)$. \emph{Смещением фракции} называем
смещение от начала до конца ее периода.
 Будем говорить, что кривая совершает \emph{поворот}  на стыке
соседних фракций одинакового порядка, если смещения этих фракций
различны.

\begin{lem}
\label{turn} Если пятерная односторонняя кривая
 делает поворот, то первая из поворотных фракций
является замедляющей, а вторая
--- ускоряющей.
\end{lem}

\begin{proof}
Подобием изменим кривую так, чтобы образы поворотных фракций
представляли собой единичные квадраты, а соответствующие им
периоды имели единичную длину. Предположим, что квадратно-линейное
отношение отмасштабированной таким образом кривой $f$ менее пяти.
  Пусть
$t_0,t_1,t_2,t_3,t_4,t_5,t_6$
--- возрастающая последовательность угловых моментов прохождения
пары поворотных фракций. Дальнейшее доказательство распадается на
два случая: поворот назад и поворот вбок.

В случае поворота назад имеем $f(t_0)=f(t_6)$. Далее лемма
\ref{2-diag}   дает $|f(t_1)-f(t_4)|^2=5$. Откуда $t_4-t_1>1$ и
$1>2-(t_4-t_1)=(t_1-t_0)+(t_6-t_4)=(t_1-t_0)+(1-(t_4-t_3))$.
Следовательно $(t_1-t_0)<(t_4-t_3)$. Тогда $t_3-t_2=t_4-t_3$ и
$t_1-t_0=t_6-t_5$ и мы получили, что первая поворотная фракция
замедляющая, а вторая --- ускоряющая.

 В
случае, бокового поворота, в силу лемм \ref{2-diag} и
\ref{str-cont} имеем $|f(t_4)-f(t_1)|^2=5$. Поэтому
$1<t_4-t_1=(t_4-t_3)+(t_3-t_2)+(t_2-t_1)$. С другой стороны
$t_3-t_1=(t_3-t_2)+(t_2-t_1)+(t_1-t_0)=1$ откуда
$t_4-t_3>t_1-t_0$.
\end{proof}
\begin{lem}
\label{turn2} Пятерная односторонняя кривая не может делать двух
поворотов подряд.
\end{lem}
\begin{proof}
Если происходит поворот на стыке первой и второй фракции, то
вторая фракция должна быть ускоряющей в силу леммы \ref{turn}, но
в этом случае в силу этой же леммы невозможен поворот на стыке
второй и третьей фракции.
\end{proof}

\begin{lem}
\label{turn+}  Если пятерная односторонняя кривая делает боковой
поворот, то фракция, предшествующая первой поворотной также как и
первая поворотная будет замедляющей.
\end{lem}
\begin{proof}
Пусть $t_0,t_1,\dots t_9$ моменты прохождения вершин трех
последовательных фракций. Если на стыке второй и третьей фракции
произошел боковой поворот, то на стыке первой и второй фракции
поворота не было, ибо два поворота подряд запрещены. Значит
$|f(t_7)-f(t_1)|^2=1^2+3^2=10$. Следовательно,
$2<t_7-t_1=(t_7-t_6)+(t_6-t_0)-(t_1-t_0)$. Поскольку $t_6-t_0=2$
отсюда следует $t_7-t_6>t_1-t_0$. Значит первая фракция ---
замедляющая.
\end{proof}

\begin{lem}
\label{1/5} Для пятерной единичной кривой расстояние между
соседними угловыми моментами лежит в диапазоне от $1/5$ до $3/5$.
\end{lem}
\begin{proof}
Предоставляется читателю.
\end{proof}

\begin{lem}
\label{turn3} Пятерная односторонняя кривая не может поворачивать
вбок после трех подряд  фракций с одинаковым смещением.
\end{lem}
\begin{proof}
Предположим противное и пусть $t_0,t_1,\dots t_{12}$
 представляет возрастающую последовательность угловых моментов соответствующих фракций.
 Тогда $|f(t_{10})-f(t_1)|^2=1^2+4^2=17$. Поэтому
$t_{10}-t_1\ge 17/5$. С другой стороны
$t_{10}-t_1=(t_{10}-t_9)+(t_9-t_0)-(t_1-t_0)$. Далее $t_9-t_0=3$,
а $t_{10}-t_9<3/5$ и $t_1-t_0>1/5$ в силу леммы \ref{1/5}. В итоге
получаем $t_{10}-t_1<3\frac25$. Противоречие.
\end{proof}

\begin{lem}
\label{turn4} Пятерная односторонняя кривая не может повернуть
назад после четырех подряд фракций с одинаковым смещением, если
первая из этих четырех фракций ускоряющая.
\end{lem}
\begin{proof}
Пусть $t_0,t_1,\dots, t_{15}$ растущая последовательность угловых
моментов рассматриваемых пяти фракций. Пятая фракция является
ускоряющей, ибо завершает поворот (лемма \ref{turn}), а первая
фракция является ускоряющей по условию. В результате
$t_{13}-t_1=4$. С другой стороны $|f(t_{13}-f(t_1)|^2=4^2+2^2=20$.
Поэтому квадратно-линейное отношение на паре $t_{13},t_1$ равно
пяти.
\end{proof}

\begin{tr}
Всякая единичная правильная односторонняя пеановская кривая имеет
квадратно-линейное отношение $\ge5$.
\end{tr}
\begin{proof}
Предположим, что  $f$ --- пятерная односторонняя кривая. Будем
считать, что $f$ имеет более четырех фракций. В случае, кривой
четвертого рода, чтобы соблюсти это условие мы просто можем
рассматривать ее как кривую шестнадцатого рода, перейдя к фракциям
второго порядка.

Пусть кривая начинает движение из левого нижнего угла, проходит
один период вверх и делает поворот направо. Так как она не может
делать двух поворотов подряд (лемма \ref{turn2}), то третий период
она продолжает движение вправо. После третьего периода она должна
повернуть вниз, иначе фракция, лежащая снизу под третьей станет
\emph{тупиком} и не сможет быть пройдена в дальнейшем, потому что
в нее можно будет попасть только справа, а выйти будет некуда. Но
на пятом периоде движение вниз уже невозможно: ниже некуда.
Приходится делать два поворота подряд, что запрещено.

Если кривая начинает движение с трех подряд периодов без
поворотов, то в дальнейшем боковые повороты становятся
невозможными в силу леммы \ref{turn3}, а повороты назад ведут в
тупик.

Пусть теперь односторонняя кривая начинает движение из левого
нижнего угла вверх и поворачивает вправо на стыке второй и третьей
фракций. Тогда третья фракция является ускоряющей в силу леммы
\ref{turn}. В четвертой фракции кривая продолжает движение вправо,
чтобы избежать двух поворотов подряд, а в стыке четвертой-пятой
фракций боковой поворот невозможен в силу леммы \ref{turn+}, а
поворот назад приводит в тупик. Поэтому движение вправо будет
продолжено в пятой фракции третий период подряд. Если в шестой
фракции продолжится движение вправо, то в дальнейшем никакие
повороты невозможны в силу лемм \ref{turn3} и \ref{turn4}.
Следовательно, в шестой фракции должен быть сделан поворот назад
(поворот вбок невозможен в силу леммы \ref{turn3}). После поворота
в шестой фракции в седьмой продолжится движение влево, которое
приводит в тупик.
\end{proof}
\begin{lem}\label{3frac}
Первые три и последние три максимальные фракции пятерной
диагональной кривой  расположены по границе ee квадрата-образа
$Q$.
\end{lem}
\begin{proof}
Вторая фракция после входа расположена на одной из сторон квадрата
$Q$, потому что есть всего четыре фракции, образы которых содержат
выход из первой фракции и только одна из них внутренняя. Но эта
фракция не может проходится второй в силу леммы \ref{str-cont}.
Далее выход из второй фракции лежит на стороне $Q$, поэтому и
третья фракция лежит на стороне квадрата $Q$.
\end{proof}
\emph{Направлением стыка} фракций называется направление от центра
первой фракции стыка на центр второй.
 \emph{ Направлением входа}
(\emph{направлением выхода}) фракции $k$-го порядка называется
направление стыка ее первых (последних) двух подфракций $k+1$-го
порядка.
\begin{lem}\label{in-out}
У пятерной диагональной  кривой  направление выхода из любой
фракции или перпендикулярно направлению входа следующей фракции
или оба эти направления параллельны пересечению образов фракций.
\end{lem}
\begin{proof}
Предположение противного приводит к тому, что тройка входных
фракций второго порядка продолжает тройку выходных в одном
направлении. Возникает шестерка фракций второго порядка вытянутая
вдоль прямой. Тогда квадрат расстояния от начала до конца этой
шестерки в шесть раз превосходит время ее прохождения.
\end{proof}

\begin{lem}\label{slow}
Диагональная пятерная единичная кривая проходит расстояние от
начала до первого внутреннего узла  не более чем за 2/5.
\end{lem}
\begin{proof}
Пусть $0=t_0,t_1,t_2,t_3=1$ суть угловые моменты первого порядка
пятерной диагональной кривой $f$. Тогда $|f(t_1)-f(t_2)|^2=2$,
поэтому $t_2-t_1>2/5$. А так как $|f(t_3)-f(t_2)|=1$, то
$t_3-t_2>1/5$. Поэтому $t_3-t_1>3/5$, откуда $t_1-t_0<2/5$
\end{proof}

 Будем
говорить, что фракция  кривой $Z$-образна, если смещение от начала
фракции до первого внутреннего угла горизонтально и будем называть
фракцию $N$-образной, если это смещение вертикально.

\begin{lem}\label{NZ}
Если пятерная диагональная кривая имеет вертикальное направление
входа, то или  первая и третья ее фракции $Z$-образны или первая
фракция замедляющая, а третья фракция ускоряющая и они имеют
разный $N-Z$-тип.
\end{lem}
\begin{proof}
Переходя, если нужно к подобной кривой, можем считать, что образы
фракций первого порядка являются единичными квадратами.
 Обозначим
через $t_i$
--- $i$-ый угловой момент первого порядка ($t_0=0$).

Если первая фракция $Z$-образна, а третья --- $N$-образна, то
$|f(t_1)-f(t_7)|^2=10$.  Пятерность кривой влечет неравенство
$t_7-t_1>2$. А так как $t_6-t_0=2$, то $t_1-t_0<t_7-t_6$, что и
означает, что первая фракция замедляющая, а третья --- ускоряющая.

Если наоборот первая фракция $N$-образна, а последняя $Z$-образна,
то $|f(t_8)-f(t_2)|^2=10$ и  $t_8-t_2>2$, откуда ввиду
$t_9-t_3=2$, получаем $t_9-t_8<t_3-t_2$, что опять-таки означает,
что первая фракция замедляющая, а третья --- ускоряющая.

Наконец, случай $N$-образных первой и третьей фракций исключен,
потому что в силу леммы \ref{slow} выполнены неравенства
$t_3-t_2<2/5$ и $t_7-t_6<2/5$. Откуда $t_7-t_3<9/5$ и
$|f(t_7)-f(t_3)|^2=10$, откуда квадратно-линейное отношение $f$ на
$t_7,t_3$ больше пяти вопреки пятерности $f$.
\end{proof}

\begin{lem}\label{4acc}
Если пятерная диагональная кривая проходит четыре фракции подряд,
двигаясь в одном и том же направлении по вертикали, то четвертая
фракция этой серии является $Z$-образной и ускоряющей.
\end{lem}
\begin{proof}
В силу леммы \ref{slow} вход, как и выход из фракции не могут
занимать более $2/5$ периода. В случае $N$-образной четвертой
секции, квадрат расстояния от начала первой секции до второй
вершины четвертой составит $17$, что сразу дает квадратно-линейное
отношение как минимум пять. Теперь сравнение второй и четвертой
секций, в силу  леммы \ref{NZ}, позволяет доказать, что четвертая
секция --- ускоряющая.
\end{proof}

\begin{tr}
Всякое фрактальная диагональная пеановская кривая имеет
квадратно-линейное отношение не менее пяти.
\end{tr}
\begin{proof}
Пусть $f$ --- пятерная диагональная кривая. Тогда
 во-первых найдется пара последовательных максимальных фракций $f$,
  такая что направление выхода из
первой фракции перпендикулярно направлению входа второй.
Действительно, предположение противного приводит в силу леммы
\ref{in-out} к заключению, что направление выхода из любой фракции
параллельно направлению входа следующей и параллельно пересечению
соседних фракций. В силу леммы \ref{3frac} найдется тройка
фракций, у которых пересечение первой и второй параллельно
пересечению второй и третьей. Поэтому для второй, а, следовательно
для любой, фракции мы получим, что направления входа и выхода
параллельны. И мы приходим с противоречием с наличием поворотов у
$f$.

Итак, изменив, если нужно, направление и ориентацию кривой, мы
получим пятерную кривую, у которой направление выхода некоторой
фракции $Q_i$ вертикально, а следующая фракция $Q_{i+1}$,
расположена над $Q_i$. В таком случае к четверке фракций второго
порядка, составленных последними тремя подфракциями $Q_i$ и первой
фракцией $Q_{i+1}$ применима лемма \ref{4acc}. То есть мы
получаем,
 что входная
подфракция $Q_{i+1}$ ускоряющая и  $Z$-образна. А так как входное
направление у $Q_{i+1}$ горизонтальное, то повернув плоскость на
прямой угол,  мы приходим к противоречию с леммой \ref{NZ}.
\end{proof}


\begin{thebibliography}{9}

\bibitem{TULA87}
Е.~В. Щепин,
\newblock Повышающие размерность отображения и непрерывная передача информации.
\newblock  {\em Вопросы чистой и прикладной математики}, том~1,
стр.   148--155. Тула. Приокское книжное издательство, 1987.



\bibitem{B} К. Е. Бауман, \newblock Коэффициент растяжения кривой
Пеано-Гильберта, препринт


\bibitem{9} E. В. Щепин, К. Е. Бауман, \newblock О кривых Пеано
фрактального рода 9. \newblock{\em Моделирование и анализ данных:
Труды факультета информационных технологий МГППУ}, (вып.1). -
Москва, РУСАВИА, 2004
\end{thebibliography}
\end{document}